# SOCP Convex Relaxation-Based Simultaneous State Estimation and Bad Data Identification

Hossein Ghassempour Aghamolki, Zhixin Miao, Lingling Fan

*Department of Electrical Engineering, University of South Florida, Tampa FL USA 33620.

Phone: 1(813)974-2031, Email: linglingfan@usf.edu.

**Abstract**

Traditional largest normalize residual (LNR) test for bad data identification relies on state estimation residuals and thus can only be implemented after running Power System State Estimation (PSSE). LNR may fail to detect bad data in leverage point measurements and multiple interacting and conforming bad data. This paper proposes an optimization problem formulation for joint state estimation and bad data identification based on second-order cone programming (SOCP) convex relaxation. $\ell_1$-norm of the sparse residuals is added in the objective function of the state estimation problem in order to recover both bad data and states simultaneously. To solve the optimization problem in polynomial time, first, SOCP convex relaxation is applied to make the problem convex. Second, least squares error (LSE)-based semidefinite programming (SDP) cutting plane method is implemented to strengthen the SOCP relaxation of PSSE. Numerical results for the proposed algorithm demonstrate simultaneous state estimation and bad data identification for networks with small sizes and thousand buses. Compared to the LNR method, the proposed method can detect bad data in leverage point measurements and multiple conforming bad data.

## 1. Introduction

Continuously monitoring the balance between power generation and consumption is the main objective for the power system operation and control. Accuracy in determination of the states of the system (i.e. voltage magnitude and voltage angle of the buses) is critical for reliable and economical operation of the power network. The results of the Power System State Estimation (PSSE) are vital for many operational discission making procedures such as unit commitment, optimal power flow (OPF), network reconfiguration, and reliability assessment [1].

To have a robust PSSE algorithm, the estimator should be able to not only provide quality estimation of the states in polynomial time, but also identify bad data measurements and remove those corrupted data to ensure the accuracy of the results.

*1.1. PSSE Solving Methods*

PSSE is a nonlinear and nonconvex problem. The nonlinearity and nonconvexity of the PSSE originate from the nonlinearity and nonconvexity of the measurement functions which relate to the power injection and line flow equations of the system. Power flow calculation is NP-hard since it needs solving a set of nonlinear polynomial quadratic equations. Also, it is known that checking the feasibility of AC power flow problem for both transmission and distribution networks is NP-hard [2, 3].

Traditionally nonlinear PSSE is formulated as a weighted least squares estimation (WLSE) and solved by iterative Gauss-Newton method [4, 1, 5]. The algorithm adapts Taylor's first-order expansion and is known to have the following issues.

- Gauss-Newton method finds a local optimum for a nonlinear WLSE problem [5]. Since nonlinear WLSE is a nonconvex optimization, it can have several local solutions. Gauss-Newton method guarantees a local solution, but not a global optimal solution.

- Sensitivity to the initialization and lack of a guaranteed convergence are the other main issues related to the Gauss-Newton method [6, 7]. It have been shown in [5] that Gauss-Newton algorithm may fail to converge due to nonlinearity or the relative large size of residual of the problem.

Hence, it is important to find more reliable solving techniques for nonlinear PSSE problems with the ability of approximating the global solution in polynomial time.

In recent years, convex optimization approach for solving nonlinear PSSE and OPF problems becomes popular. Two well-known convex relaxation methods are semidefinite programming (SDP) and second order cone programming (SOCP). SDP relaxation is first applied in OPF in [8]. Lavaei and Low have shown that SDP relaxation of OPF can provide a tight bound for OPF problems [9]. A summary report on the recent advances in convex relaxations of the OPF problem can be found in [10] and [11].

Using SDP relaxation for PSSE problems has been seen in the literature [12, 7, 13, 14, 15, 16]. The disadvantage of SDP relaxation for PSSE problem is that it has scalability issue as demonstrated with computing experiments in [17]. For large-size power networks, computing time of the SDP relaxation is significant and decomposition techniques become an option [18]. The power networks used in [12], [7], [13], [14], and [15] are limited to a 14-bus system, a 30-bus system, a 6-bus system, a 6-bus system, and a 30-bus system, respectively.

To be able to handle large-size networks, distributed computing techniques are adopted for SDP relaxation of PSSE. For example, [19, 20, 21] adopted distributed optimization for SDP relaxation of PSSE. The main idea is to decompose a large-scale power network into subsystems. For each subsystem, an optimization problem will be solved. Among the subsystems, information will be exchanged and an iterative procedure will be carried out until convergence. For example, Lagrangian method is employed in [19] and alternating direction method of multipliers (ADMM) is adopted in [20, 21]. Distributed algorithms have also been applied to Gauss-Newton method, see [22, 23, 24].

Besides distributed computing, more computationally efficient convex relaxations such as linear programming relaxation [25], quadratic programming relaxation [26], and SOCP relaxation [17] have the potentials to be applicable in PSSE. Among those methods, SOCP has demonstrated a great potential to achieve a tighter gap with more computational efficiency [17, 27].

Using SOCP relaxation for solving power flow problems was first initiated in [28] and [29]. In [6], SOCP relaxation for PSSE is investigated. To have the relaxation exact, the objective function of a WLSE problem is modified to include voltage phasor information of a spanning tree. Jabr and Pal used the technique of variable change similar to that in SOCP formulation to linearize the power flow equations in PSSE [30]. SOCP relaxation was not adopted in [30] since SOCP relaxation can be weak for meshed networks [27, 31]. Therefore, a nonconvex optimization problem is formulated in [30] and interior point method (IPM) solver is adopted to find a local optimal solution.

Recently, cycle based constraints have been introduced to enhance SOCP relaxation [27, 31]. In particular, [27] has suggested three approaches for strengthening SOCP relaxation for OPF problems using cycle based feasibility constraints. The SDP separation approach is claimed to provide the best strengthening effect. The main idea is to add linear inequality constraints to reduce the feasibility region.

In this paper we adopt a new strengthened SOCP relaxation method [32] for PSSE problems. In our method, affine inequality constraints instead of linear inequality constraints are used. Our strengthening method provides a similar strengthening effect. Detailed mathematic derivation and theoretic investigation on this algorithm can be found in the authors' paper [32], where the algorithm was applied for OPF.

*1.2. Bad data Detection*

Corrupted data usually exist in power system measurements due to device failures, communication system problems or cyber attacks. PSSE is expected to detect and get rid of bad data to guarantee the accuracy of the estimation results. Bad data detection and identification have been widely reported in literature. Schweppe in [33, 34] has introduced PSSE and bad data detection with WLSE formulation. Different methods such as weighted and normalized residuals, sum of squared residuals, and non-quadratic criteria for solving PSSE and bad data identification problems are introduced in [35]. In the literature, $\chi^2$ test was introduced to show if there are bad data while Largest Normalized Residual (LNR) test shows exactly which measurement contains bad data [35, 1].

There are three issues related to LNR tests.

- LNR test relies on the state estimation residuals and thus can only be implemented after running PSSE. In case of detecting any bad datum, PSSE has to run again after discarding that bad datum.

- LNR test may fail to detect multiple interacting conforming bad data [1].

- If leverage measurements contain bad data, the LNR algorithm would not be able to detect it [1, 36, 37].

In order to overcome the limitations of LNR in the detection of multiple interacting conforming bad data, several research works have been carried out in the literature. In [38], hypothesis testing was introduced to identify multiple interacting bad data while decomposing scheme with least median squares estimation was suggested in [39]. More recently, [40] has suggested a modified LNR algorithm to specifically deal with multiple interacting bad data identification.

In the context of leverage points bad data detection, projection statistics algorithm has been introduced in [1] and [36]. Most recently Pal and Majumdar introduced diagnostic robust generalized potential algorithm in [37] for separating leverage measurements in power system.

The aforementioned research all rely on separated state estimation and bad data identification. Recently proposed robust state estimation methods [41, 7, 42] can conduct simultaneous state estimation and bad data identification through joint optimization problem formulation. In fact, the identification problem for both state variables and bad data (of the same dimension as the measurement vector) is an under-determined problem [7]. If the bad data are sparse signals, then the estimation problem becomes over-determined.

The application of sparse signal recovery for bad data identification gains lots of attentions in recent years. *Sparse Residual Estimator* (SRE) with $\ell_1$-norm optimization was first introduced in [43] for recovering sparse signals from unreliable sensors in the network by using sparse matrix characteristics to create an M-Huber estimator. In [41], SRE application has been expanded by using the method for developing joint *linear* state estimation and bad data identification. In [7], SDP relaxation of the robust nonlinear state estimation problem is solved and the power network is limited to 30 buses. In [42], Gauss-Newton type of iterative solver is employed to find the local solution of the robust state estimation problem. The test system is also limited to 30 buses.

The objective of this paper is to formulate a robust nonlinear PSSE problem to simultaneously identify state variables and bad data and further solve the problem for large-scale networks with approximate global solutions using strengthened SOCP relaxation.

*1.3. Contributions*

Motivated by [41, 7] and our previous research in strengthened SOCP relaxation based OPF [32], we contrive convex optimization framework for PSSE problems using SOCP relaxation with cycle based constraints for the first time. Relating the feasibility constraints in SOCP relaxation of PSSE to every cycle in a cycle basis, makes it possible to generate valid affine inequalities (or cuts) to reduce the search space. The advantage of the proposed method for PSSE is demonstrated by comparing the suggested method with the traditional iterative Gauss-Newton algorithm.

Sparse residual estimator (SRE) has been used for introducing joint state estimation and bad data identification only for linear system in literature [43, 41]. In fact [42] showed that, in case of nonlinear and nonconvex measurement functions, there is no guarantee that SRE algorithm would recover sparse error. Motivated by those research, joint PSSE and bad data identification algorithm is proposed in this paper by using SRE and strengthened SOCP relaxation formulation. The advantage of proposed method then compared with traditional WLSE-based PSSE and traditional LNR bad data detection. Two example test studies have been carried out to show the effectiveness of proposed joint optimization algorithm in the detection of leverage points bad data and multiple conforming bad data. Finally, the robustness of the method is tested against 13 NESTA V6.0 test cases from a 3-bus system to a 1354-bus systems [44].

The rest of the paper is organized as follows. Section 2 presents state estimation formulation and introduces joint state estimation and bad data identification as a constrained optimization problem. SOCP relaxation of the problem is also presented. Section 3 presents the LSE based SDP cutting plane method that strengthens the SOCP relaxation. In Section 4, case studies are carried out to test and verify the effectiveness and robustness of the proposed algorithm. Finally Section 5 concludes the paper.

## 2. State estimation and the proposed co-optimization problem formulation

This section presents state estimation in its standard formulation with two types of objective functions. SOCP relaxation is introduced, followed by joint state estimation and bad data identification as a constrained optimization problem.

Consider a power network with $\mathcal{N} = (\mathcal{B}, \mathcal{L})$ where $\mathcal{B}$ is the set of buses of the system and $\mathcal{L}$ is the set of transmission lines. Assume that $\mathcal{G}$ is the set of buses on the network with generators connected to them and $Y$ is the admittance matrix of the system with $Y_{ij} = G_{ij} + jB_{ij}$ as a complex rectangular representation of its elements. For a transmission line connection bus $i$ and bus $j$, $y_{ij} = g_{ij} + jb_{ij}$ represents its series admittance in the rectangular form and $b_{sh}$ represents the charging susceptance in its $\pi$ equivalent model.

*2.1. State estimation standard formulation*

Voltage magnitudes ($|V|$) and angles ($\delta$) of every bus of the power network are usually assumed as state variables ($x$) in PSSE. PSSE uses measurement data from a snapshot of power system operation to produce a reliable estimation of the transmission line flow and voltage of the buses. There are two different methods widely used for PSSE: *Weighted Least Square Error (WLSE) Estimation* and *Weighted Least Absolute Value (WLAV) estimation*. WLSE minimizes the sum (or weighted sum) of the squares of the errors between the measured values and the measurement function (see (1)). Consequently, the cost function of the PSSE problem can be represented by the $\ell_2$-norm of the errors between the measurement functions values and the measurements data. Thus, WLSE cost function is a quadratic convex function.

$$\text{minimize} \quad \sum_{i=1}^{m} \left(\frac{r_i}{\sigma_i}\right)^2 \tag{1a}$$

$$\text{subject to} \quad z_i = h_i(x) + r_i, i = 1, \cdots, m \tag{1b}$$

where $x \in \mathbb{R}^{2 \times |\mathcal{B}|}$, $z \in \mathbb{R}^m$, $z_i$ is the $i$th measurement, $h_i(x)$ is the nonlinear function relating the state vector $x$ to $z_i$, $r_i$ represents the residual, $\sigma_i$ represents $i$th measurement's deviation.

Compared to WLSE, WLAV estimation minimizes the sum of the weighted absolute errors between measured values and measurement function outputs (see (2)). The objective function of WLAV is not differentiable.

$$\text{minimize} \quad \sum_{i=1}^{m} \left|\frac{r_i}{\sigma_i}\right| \tag{2a}$$

$$\text{subject to} \quad z_i = h_i(x) + r_i, i = 1, \cdots, m \tag{2b}$$

There are additional constraints, *e.g.*, zero-injection pseudo-measurement constraints (3a) and/or inequality constraints (3b) (e.g., constraints related to the direction of the power injection or flow which is only used in conjunction with the current magnitude measurements).

$$r_i = 0, \quad i \in \mathcal{P} \tag{3a}$$

$$z_i^{\min} \leq h_i(x) \leq z_i^{\max} \tag{3b}$$

where $\mathcal{P}$ notates the set of the pseudo-measurement meters.

Measurements of the system usually consist of voltage magnitudes of the buses, real and reactive power injections, and real and reactive power flow of transmission lines. Since AC power flow is a nonlinear, non-convex problem, PSSE will be nonlinear and non-convex optimization as well. The power injection equations at Bus $i$ can be written as (4a)-(4b).

$$P_i^g - P_i^d = \sum_{j=1}^{N} |V_i||V_j|[G_{ij}\cos(\delta_{ij}) + B_{ij}\sin(\delta_{ij})] \tag{4a}$$

$$Q_i^g - Q_i^d = \sum_{j=1}^{N} |V_i||V_j|[G_{ij}\sin(\delta_{ij}) - B_{ij}\cos(\delta_{ij})] \tag{4b}$$

The line flow equations at a branch connecting Bus $i$ and Bus $j$ are written as (5a)-(5b).

$$P_{ij} = g_{ij}|V_i|^2 - |V_i||V_j|[g_{ij}\cos(\delta_{ij}) + b_{ij}\sin(\delta_{ij})] \tag{5a}$$

$$Q_{ij} = -\left(b_{ij} + \frac{b_{sh}}{2}\right)|V_i|^2 + |V_i||V_j|[b_{ij}\cos(\delta_{ij}) - g_{ij}\sin(\delta_{ij})] \tag{5b}$$

Superscripts $^g$ and $^d$ notate generation and load consumption respectively, and $\delta_{ij} = \delta_i - \delta_j$.

The state estimation problem will give estimate of the state: $|V|$ vector and $\delta$ vector based on the measurement data consisting of voltage magnitude measurements, power injection measurements, and line flow meaurements.

2.2. SOCP relaxation of the state estimation

From equations (4a)-(5b), it can be observed that the polar form of voltage phasor expression leads to nonlinearity and non-convexity of the measurement functions. The voltage phasor will be expressed in the reactangular form:

$$\overline{V}_i = e_i + jf_i.$$

The nonlinear terms can be converted to linear convex term by defining new variables $T$ and $S$ for all lines of the power network.

$$T_{ii} = e_i^2 + f_i^2 = |V_i|^2, i \in \mathcal{B} \tag{6a}$$

$$T_{ij} = e_i e_j + f_i f_j = |V_i||V_j|\cos(\delta_i - \delta_j), (i,j) \in \mathcal{L} \tag{6b}$$

$$S_{ij} = e_i f_j - f_i e_j = -|V_i||V_j|\sin(\delta_i - \delta_j), (i,j) \in \mathcal{L} \tag{6c}$$

Substituting the new variables into the equations (4a)-(5b), the power injection equations can be represented as follows.

$$P_i^g - P_i^d = G_{ii} T_{ii} + \sum_{j \in a_i} [G_{ij} T_{ij} - B_{ij} S_{ij}] \tag{7a}$$

$$Q_i^g - Q_i^d = -B_{ii} T_{ii} - \sum_{j \in a_i} [B_{ij} T_{ij} + G_{ij} S_{ij}] \tag{7b}$$

where $a_i$ is the set of the adjacent buses that have branches connected to Bus $i$. The line flow expressions at branch $k$ connecting bus $i$ and bus $j$ are as follows.

$$P_{ij} = g_{ij} T_{ii} - g_{ij} T_{ij} + b_{ij} S_{ij} \tag{8a}$$

$$Q_{ij} = -\left(b_{ij} + \frac{b_{sh}}{2}\right) T_{ii} + b_{ij} T_{ij} + g_{ij} S_{ij} \tag{8b}$$

The following relations between new introduced variables have to be held:

$$T_{ij} = T_{ji}, \quad S_{ij} = -S_{ji} \tag{9a}$$

$$T_{ij}^2 + S_{ij}^2 = T_{ii} T_{jj} \tag{9b}$$

The last equation in (9a) is a surface of a cone and is a non-convex constraint. We can apply relaxation by changing this constraint from equality to inequality.

$$T_{ij}^2 + S_{ij}^2 \leqslant T_{ii} T_{jj} \tag{10}$$

The above inequality describes a second-order cone, hence the name SOCP relaxation. The above relaxation was first introduced in [45] and [28]. The latter showed that the SOCP relaxation is exact for radial networks.

The SOCP relaxation of PSSE will give estimate of $S$ and $T$ based on the measurement data. A slight change is assumed for voltage magnitudes. Instead of directly using a voltage magnitude in the $z$ vector, in the SOCP relaxation problem, we use its square in the measurement data. This gives us a linear measurement function:

$$(|V|_i^{\text{meas}})^2 = T_{ii} + r,$$

where superscript $^{\text{meas}}$ notates measurement.

In the conventional state estimation problems, the voltage magnitude measurement's standard deviation is given based on the meter accuracy. Assume that $|V|_i^{\text{meas}} = |V|^{\text{true}} + e$, where the error $e$ is Gaussian noise with zero mean and variance as $\sigma^2$: $e \sim N(0, \sigma^2)$.

Then

$$(|V|_i^{\text{meas}})^2 \approx (|V|^{\text{true}})^2 + 2(|V|^{\text{true}})e.$$

A system's voltage is close to 1 pu. Hence we may assume that $r$ is also a Guassian distribution with its deviation as $2\sigma$.

The optimal solution for this problem may not be feasible and the relaxation is not exact due to two reasons.

- Constraint (10) is not binding, i.e., : $T_{ij}^2 + S_{ij}^2 < T_{ii}T_{jj}$.

- For meshed networks with cycles, Kirchhoff's Voltage Law (KVL) of a loop is not satisfied. This can be manifested as the sum of the angle differences across a cycle is not zero [27]: $\sum_{(i,j) \in C} \tan^{-1}\left(\frac{-S_{ij}}{T_{ij}}\right) \neq 0$.

To take care of Kirchhoffs Voltage Law and make the relaxation exact for the meshed network, [29] proposed imposing following constraints to the SOCP formulation:

$$\delta_j - \delta_i = \tan^{-1}\left(\frac{S_{ij}}{T_{ij}}\right). \tag{11}$$

Note that, (11) contains tangent function which makes it a non-convex constraint. Therefore, the constraint has to be linearized in order to achieve convex SOCP relaxation of the alternative formulation. A heuristic iterative procedure is introduced in [29] to incorporate the linearized angle constraint. However, the introduced linear constraints may exclude feasible solutions and the final solution is only an approximation. In Section III, we will show alternative ways to handle the issue related to meshed networks.

*2.3. Proposed joint state estimation and bad data identification problem formulation*

Nonlinear measurement function can be represented by the following equation:

$$z = h(x) + w \tag{12}$$

where $z$ is the measurement vector, $x$ is the state variables vector, $h(x)$ is the measurement function coefficient and $w \in \mathbb{R}^m, w_i \sim N(0,1)$ represents Gaussian noises of standard distribution. If corrupted data exist in the measurements, a sparse vector $o$ can be added as an unknown vector which only has non-zero elements $o(i)$ if $z(i)$ contains bad data [42, 43, 41]. In this case, the new measurement model can be represented as (13).

$$z = h(x) + o + w \tag{13}$$

Joint estimation of $x$ and $o$ can reveal states while identifying the corrupted data. [42] shows that relying on the sparsity of $o$, if a list of $\tau_0$ faulty measurements are expected, ideally a combination of $\ell_0$-pseudonorm and $\ell_2$-norm as shown in (14) could successfully recover $x$ and $o$.

$$\begin{aligned} \min_{x \in X, o} & \ ||z - h(x) - o||_2^2 \\ \text{s.t.} & \ ||o||_0 \leq \tau_0 \end{aligned} \tag{14}$$

In the above equation, $\tau_0$ is the number of corrupted measurements. The problem is $\ell_0$-pseudonorm in (14), which renders NP-hard and makes it computationally impossible to solve the optimization problem for large scale systems. To make the problem computationally efficient, a well-known convex $\ell_1$-norm relaxation can be applied to above constraint [42], [41].

$$\begin{aligned} \min_{x \in X, o} & \ ||z - h(x) - o||_2^2 \\ \text{s.t.} & \ ||o||_1 \leq \tau_1 \end{aligned} \tag{15}$$

where $\tau_1$ is a positive number and depends to the number of corrupted measurements and the probability distribution function of error in those measurements. The Lagrangian relaxation of the constraint leads to the co-optimization represented in (16), which allows joint state estimation and bad data identification simultaneously.

$$(\hat{x}, \hat{o}) \in \arg\min_{x \in X, o} \left\{ ||z - h(x) - o||_2^2 + \lambda ||o||_1 \right\} \tag{16}$$

where $\lambda$ is a positive parameter. Based on the Lagrangian formulation in (16), it can be anticipated that if $\lambda \to \infty$, then $\hat{o}$ becomes zero and the joint optimization reduces to WLSE. On the contrary, it has been shown in [46] and [47] that when $\lambda \to 0^+$, the solution of the joint formulation coincides with WLAV. Also it has been shown in the literature that for a finite $\lambda > 0$, joint formulation will respond to Huber's M-estimator; see [41] and [48] for more details. With the assumption that $w$ relates to Gaussian noise, choosing $\lambda = 1.34$ makes the estimator 95% asymptotically efficient for uncorrupted measurements [4]. The sensitivity of the joint formulation to $\lambda$ is further discussed by carrying out some examples in the case study section.

The above joint optimization is convex if the measurement functions $h(.)$ are convex. In standard AC state estimation, measurement functions are nonlinear and non-convex as shown in equations (4a)-(5b). [42] showed that, in case of having nonlinear and nonconvex measurement functions, there is no guarantee for the joint optimization algorithm to recover sparse error. Therefore, for the first time, in this paper we have adopted SOCP relaxation to build a convex optimization problem with an $\ell_1$ norm constraint to make it possible for the joint optimization formulation to be applied for nonlinear PSSE. Also, the use of our new LSE-based SDP cuts strengthens SOCP relaxation, avoids the need for a local solver or iterative linearization in joint optimization solving, and leads to a more accurate estimation and identification result. Note that the state variables in our constrained optimization problem are no longer voltage magnitudes and angles, instead, $S$ and $T$.

## 2.4. Feasibility check and voltage phasor recovery

Two steps are followed after $S$ and $T$ are obtained from the strengthened SOCP relaxation of state estimation problem, with Step 1 to check if the solution is a feasible solution and Step 2 to recover the voltage phasor from $S$ and $T$ vectors.

In Step 1, two feasibility indices $\text{err}_{\text{socp}}$ and $\text{err}_{\text{cycle}}$ are used for feasibility check. The first index checks if (9b) is satisfied while the second index checks if the sum of the phase angle differences across every cycle in the grid is zero.

An error `err_branch_k` $(= T_{ii}T_{jj} - T_{ij}^2 - S_{ij}^2)$ is defined for each branch to notate the difference between $T_{ii}T_{jj}$ and $T_{ij}^2 + S_{ij}^2$. Define the sum of the errors as

$$\text{err}_{\text{socp}} = \sum_k \texttt{err\_branch\_k}.$$

An error `err_cycle_k` $\left(=\left|\sum_{(i,j)\in C}\delta_{ij}\right|\right)$ is defined for each cycle in the cycle basis of the power network.

$$\text{err}_{\text{cycle}} = \sum_k \texttt{err\_cycle\_k}.$$

If both indices are zero, the solution from the SOCP relaxation is a feasible solution.

Whether the solution is a feasible or not, we recover the voltage phasors. An infeasible solution can be used as the initial start for the Gauss-Newton method and a more accurate estimation can be obtained according to [12].

To recover $|V|$ and $\delta$ from the solution of the SOCP relaxation problem ($S$ and $T$), the following recovery steps will be adopted. For Bus $i$, the voltage magnitude $|V_i|$ can be computed from $T_{ii}$:

$$|V_i| = \sqrt{T_{ii}}.$$

For branch $k$ connecting bus $i$ and bus $j$, based on $S_{ij}$ and $T_{ij}$, we are able to find the difference between these two voltage phasor angles:

$$\delta_j - \delta_i = \tan^{-1}\left(\frac{S_{ij}}{T_{ij}}\right)$$

Aggregating the equations related to all branches, we should arrive at a matrix/vector expression that relates the phase angle vector to the vector $b \in \mathbb{R}^m$, where $m$ is the total number of branches in the power grid and $b_k = \tan^{-1}\left(\frac{S_{ij}}{T_{ij}}\right)$:

$$M\delta = b$$

where $M \in \mathbb{R}^{|\mathcal{L}|\times|\mathcal{B}|}$ is the branch to bus incidence matrix. For the $k$-th row of $M$, all components except those in the $j$th column and $i$th column are not zero.

$$M_{ki} = -1, \quad M_{kj} = 1$$

if Branch $k$ is connecting Bus $i$ to Bus $j$.

$M$ is generally not a square matrix. Using MATLAB function `linsolve(M, b)`, we may find the phase angle vector $\delta$. Further, if the reference bus phase angle is not zero, this vector will be updated to have the reference bus phase angle zero while the rest of the phase angles are all subject to a shift of $-\delta_{\text{ref}}$, where ref is the reference bus index.

## 3. Strengthening SOCP relaxation using LSE-based SDP cuts

SOCP relaxation may be weak for meshed power networks due to its ignoring the feasibility constraints in (11). On the other hand, SDP relaxation is known to produce very tight lower bound and providing quality guarantees for power flow optimization. However, SDP algorithms have scalability issues [31]. In order to have the benefit of efficient computing of SOCP relaxation together with the benefit of tight lower bound results of SDP relaxation, SDP cuts are introduced to strengthen SOCP relaxation.

The idea is to strength SOCP relaxation by adding inequality constraints to separate SOCP solutions from the feasible region of the SDP relaxation. This separation idea is first proposed in [27]. Linear inequality constraints are generated using the method proposed in [27]. In this paper, we use a method proposed in [32] to obtain affine inequality constraints using LSE-based method.

### 3.1. SDP relaxation

In SDP relaxation, rectangular expressions are used to represent the voltage phasors.

$$\overline{V}_i = |V_i|\angle\delta_i = |V_i|\cos\delta_i + j|V_i|\sin\delta_i = e_i + jf_i \tag{17}$$

A matrix $W$ is defined as follows

$$W = \begin{bmatrix} e \\ f \end{bmatrix} \begin{bmatrix} e^T & f^T \end{bmatrix} \tag{18}$$

where $f = (f_1, f_2, \cdots, f_n)^T$, $e = (e_1, e_2, \cdots, e_n)^T$.

It is obvious to find the following characteristics:

$$W = W^T \text{ and } W \succeq 0 \tag{19}$$

$W \succeq 0$ means that this matrix is positive semi-definite.

The power injection constraints will be shown to be linear with the elements of $W$. Define

$$i' = i + |\mathcal{B}|, \quad j' = j + |\mathcal{B}|, \tag{20}$$

where $|.|$ notates the cardinality of a set and $\mathcal{B}$ is the number of buses in the grid.

$$P_i^g - P_i^d = \sum_{j=1}^{N} G_{ij}(W_{ij} + W_{i'j'}) + B_{ij}(W_{ji'} - W_{ij'})$$

$$Q_i^g - Q_i^d = \sum_{j=1}^{N} (G_{ij}(W_{i'j} - W_{ij'}) + B_{ij}(W_{ij} + W_{i'j'}))$$

The above expressions indicate that the equality constraints of power injection are linear in terms of $W$. Without the rank 1 constraint rank($W$) = 1, this problem is a convex problem and a semidefinite programming (SDP) problem. This relaxation is named as SDP relaxation [8].

3.2. SDP cuts

The variables utilized in SOCP relaxations are $S$ and $T$ and they have the following relations with $W$.

$$T_{ij} = e_i e_j + f_i f_j = W_{ij} + W_{i'j'} \qquad (i,j) \in \mathcal{L} \qquad (21a)$$

$$S_{ij} = e_i f_j - e_j f_i = W_{ij'} - W_{ji'} \qquad (i,j) \in \mathcal{L} \qquad (21b)$$

$$T_{ii} = e_i^2 + f_i^2 = W_{ii} + W_{i'i'} \qquad i \in \mathcal{B} \qquad (21c)$$

Based on above equations, For every $T_{ij}$, $S_{ij}$ and $T_{ii}$, we can express them to be the Frobenius product related to the matrix $W$ by the following equation:

$$T_{ij} = A_l \bullet W = \text{Trace}(A_l W^T) \qquad (22)$$

where $\bullet$ denotes Frobenious product.

Therefore, for any solution $(T^*, S^*)$ of SOCP relaxation, if we could find a positive semidefinite $W^*$ that satisfies equations in (21a)-(21c), then $(T^*, S^*)$ belong to the SDP feasible region and it is an optimal solution for SDP relaxation. On the contrary, if such $W^*$ does not exist, the optimal solution $(T^*, S^*)$ of SOCP relaxation does not belong to the SDP feasible region and thus, it is not a feasible solution for PSSE problem. Therefore, we can separate $(T^*, S^*)$ from SDP feasible region by adding a set of inequality constraints to the SOCP relaxation.

A challenge is that finding a full size $W$ needs solving an SDP optimization with the matrix of the same size of the original SDP relaxation problem. This is a very time consuming task. Consequently, instead of full size matrix search, we can separate the solution to the SOCP relaxation problem over every cycle in a cycle basis. For any chordless cycle $C$ belongs to the cycle basis, we are looking to find corresponding submatrix $\tilde{W}$ of $W$. In this way, the separation will be very efficient and effective. The SDP feasible region we are looking for can be defined as follows.

$$\left\{ \begin{array}{c} \mathcal{S} := z \in \mathbb{R}^{3|C|} : \exists \tilde{W} \in \mathbb{R}^{2|C| \times 2|C|} \text{ s.t.} \\ -z_l + A_l \bullet \tilde{W} = 0 \quad \forall l \in \mathcal{L}, \quad \tilde{W} \succeq 0 \end{array} \right\} \tag{23}$$

To find proper inequality constraints to separate SOCP relaxation results from the SDP's feasible region, we recently introduced a new LSE-based method. The philosophy is explained by Fig. 1.

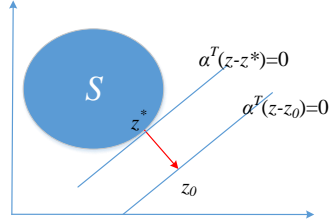

Figure 1: LSE based SDP cut will add $\alpha^T(z - z^*) \leq 0$ as cuts to the SOCP problem, where $\alpha^T = z_0 - z^*$.

For any optimal solution of SOCP relaxation $z_0$ that belongs to cycle $C$ in a cycle basis, first, we will find the shortest distance from $z_0$ to the set $S$ where $z^*$, and $\tilde{W}^*$ are the corresponding variables found in $S$. Therefore, a small-scale LSE optimization problem over a cycle will give us the corresponding $z^*$, and $\tilde{W}^*$ for set $S$.

$$\min_{z, \tilde{W}} \quad \|z_0 - z\|^2 \tag{24a}$$

$$st. \quad z_l = \text{Trace}(A_l \tilde{W}^T), l \in L \tag{24b}$$

$$\tilde{W} \succeq 0 \tag{24c}$$

where $\tilde{W}$ is a corresponding submatrix of $W$ for Cycle $C$, $L$ is the set of all indices of $z$. If the norm of $z_0 - z^*$ is zero, that means $z_0$ belongs to the SDP set and $z_0$ meets the requirement of cycle constraint. Therefore, $\alpha = z_0 - z^* = 0$ and no cuts will be added to the SOCP relaxation problem. On the other hand, if $z_0$ does not belong to the SDP feasible region, then $\alpha \not\equiv 0$ and $\alpha^T(z - z^*) \leq 0$ inequality constraints will be added to the original SOCP problem. This procedure is applied for every cycle in a cycle basis. Therefore, in each iteration several cuts will add to the original problem. Our studies show that even 3 to 4 iterations are enough to reach to a quality result for PSSE with SOCP relaxation.

The mathematical details of the SDP cut-based strengthened SOCP were documented in [32]. Comparison of this approach and the strengthening method in [27] can also be found in [32].

## 4. Case Study and Numerical results

In this section, case studies are presented to show the effectiveness of the proposed simultaneous state estimation and bad data identification algorithm. The algorithm is programmed and implemented using CVX toolbox of the MATLAB [49]. MATLAB have been running on a Core2Duo, 3.00 GHz PC with 6.00 GB of RAM. MOSEK solver [50] is used for solving convex optimization. In each scenario of case studies, power flow results from MATPOWER [51] have been used as the true values of the measurements. Also, MATPOWER's Gauss-Newton state estimation function `run_se` has been modified to provide comparison. Noises are represented by random Gaussian distribution values with zero mean and standard deviation of 0.01 pu for power and 0.005 pu for voltage.

### 4.1. SOCP versus Strengthen SOCP for WLSE

The first test is conducted on IEEE 14-bus system for WLSE only. This system has 7 cycles in its cycle basis. The objective of this test is to demonstrate the advantage of strengthened SOCP over SOCP relaxation in providing a solution much closer to the feasible region. SOCP relaxation of WLSE is first solved. The solution is used to create 7 SDP cuts. The SOCP problem is now modified to include the 7 SDP cuts. This is the first iteration. Once we obtain the solution, we may create another 7 cuts and solve the strengthened SOCP with 14 cuts.

The feasibility indices from SOCP relaxation, SOCP relaxation with SDP cuts for one iteration, and SOCP relaxation with SDP cuts for two iterations are compared for 50 Monte Carlo runs of the WLSE problem in Fig. 2. At each run, the measurement data were created by adding Gaussian noise on top of the true values.

It can be seen clearly that the strengthened SOCP relaxation provides a solution much closer to the feasible region.

### 4.2. Performance of the strengthened SOCP on WLSE-based PSSE

The proposed estimator based on strengthened SOCP relaxation with LSE-SDP cuts are tested on IEEE-14, 30 and 118 bus test systems. The results of the proposed estimator are compared with those obtained

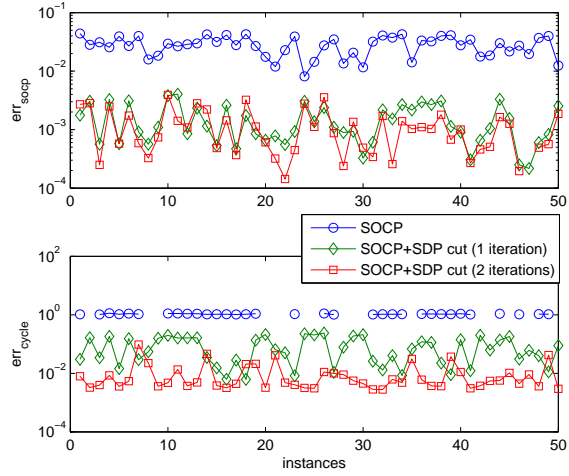

Figure 2: The feasibility indices of the SOCP relaxation, SOCP relaxation with 7 SDP cuts and SOCP relaxation with 14 SDP cuts of the WLSE for IEEE 14-bus system.

from traditional Gauss-Newton method. Equations (25a)-(25b) shows the performance indices.

$$\text{RMS--VE} = \sqrt{\frac{\sum_{i=1}^{N}(V_i^{\text{true}} - V_i^e)^2}{N}} \tag{25a}$$

$$\text{RMS--AE} = \sqrt{\frac{\sum_{i=1}^{N}(\delta_i^{\text{true}} - \delta_i^e)^2}{N}} \tag{25b}$$

where superscript $^{\text{true}}$ stands for true value and $e$ stands for estimation. $N$ is the number of buses in the network.

In order to ensure the observability of the system, the conventional measurement sets have been adopted from [30]. For the IEEE-118 bus system, the measurement set consists of (i) the real and reactive power flow for each end of transmission lines, (ii) the real and reactive power injection at every bus, and (iii) voltage magnitude at every bus. The performance indices are the Mean Square Error (MSE), the maximum of root mean square (RMS) voltage error (RMS-VE), and the maximum of RMS angle error (RMS-AE).

The results are presented in Figs. 3-5. For every system, 50 instances of measurement data were created. The two algorithms were tested for the same measurement data at each run. It can be seen from the three figures that when the size of the network becomes larger, strengthened SOCP (with 3 iterations) shows better performance compared to the Gauss-Newton method.

Since the algorithm needs to be solved iteratively, its computational time is slightly higher than the Gauss-Newton method. The computational time of the two methods we tested are compared in Table 1.

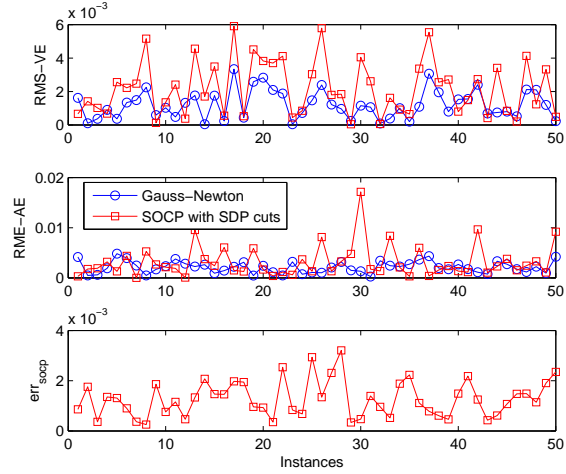

Figure 3: Comparison of the proposed method and the Gauss-Newton method for the IEEE-14 bus system.

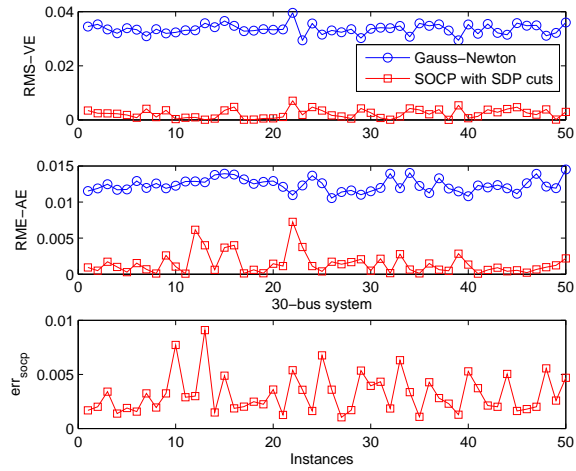

Figure 4: Comparison of the proposed method and the Gauss-Newton method for the IEEE-30 bus system.
18

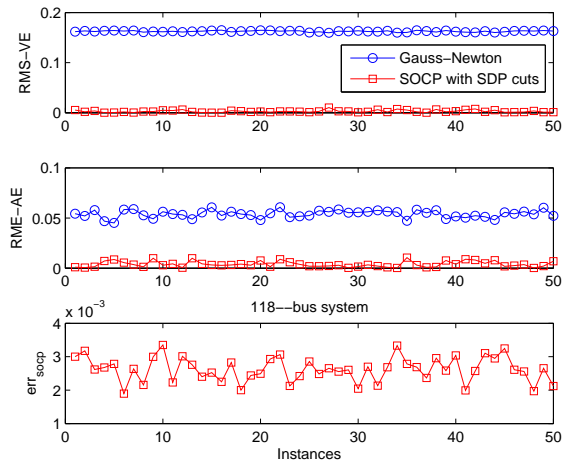

Figure 5: Comparison of the proposed method and the Gauss-Newton method for the IEEE-118 bus system.

The computational time is broke down for the strengthened SOCP method to have the computing time for solving the first SOCP relaxation problem, creating SDP cuts, and solving the strengthened SOCP relaxation problem with the SDP cuts. Only one iteration is counted in this table for the total computing time. In the step of SDP cuts creation, for each cycle in the cycle basis of the system, an SDP problem in the form (24) is solved. The three networks have 7, 12, and 62 cycles in their cycle basis, respectively. The more cycles, the more computing time is needed.

Table 1: Computational time in seconds comparison between two methods

| Cases | PSSE-SDP Cuts | | | | Newton | #Cycles |
|---|---|---|---|---|---|---|
| | SOCP | SDP cuts | SOCP-SDP cuts | total | | |
| IEEE 14 | 0.6839 | 0.8566 | 0.6930 | 2.2335 | 0.08 | 7 |
| IEEE 30 | 1.3251 | 1.5862 | 1.2556 | 4.1669 | 0.1172 | 12 |
| IEEE 118 | 7.1309 | 9.5163 | 7.1421 | 23.7893 | 1.2043 | 62 |

### 4.3. Performance of the joint optimization problem formulation against multiple conforming bad data

This section investigates the effectiveness of the proposed joint optimization method in the presence of interacting and conforming multiple bad data. In order to do that, Example 5.5 of the reference [1] have been implemented on NESTA [44] three-bus system to simulate multiple interacting conforming bad data.

A brief definition of interacting conforming bad data [1] is also explained here. The measurement residual vector $r$ can be expressed by the product of the sensitivity matrix $S$ and the error vector $e$:

Consider the linearized measurement equations:

$$\Delta z = H \Delta x + e,$$

where $e$ is assumed to have normal distribution: $E(e) = 0$, and $cov(e) = R$. Then the weighted least square (WLSE) estimator is:

$$\Delta \hat{x} = (H^T R^{-1} H)^{-1} H^T R^{-1} \Delta z = G^{-1} H^T R^{-1} \Delta z.$$

The estimated value of $\Delta z$ is:

$$\Delta \hat{z} = H \Delta \hat{x} = K \Delta z,$$

where $K = H G^{-1} H^T R^{-1}$.

The residual vector $r$ can then be expressed by the error vector.

$$\begin{aligned} r = \Delta z - \Delta \hat{z} &= (I - K) \Delta z = (I - K)(H \Delta x + e) \\ &= \underbrace{(I - K)}_{S} e \quad \text{Since } KH = H \end{aligned} \tag{26}$$

If $S_{ij}$ is very small compared to $S_{ii}$ and $S_{jj}$, then measurements $i$ and $j$ are noninteracting. That is, error in one measurement will not influence the estimate of the other variable. If $S_{ij}$ is large, then the two measurements are interacting. The interacting measurements are then classified as nonconforming and conforming. If the errors in the two measurements are not consistent with each other, then the LNR test may still identify the bad data. This type is interacting nonconforming. If the errors in the two measurements are in agreement, then the largest normalized residual test will fail. This type is interacting conforming.

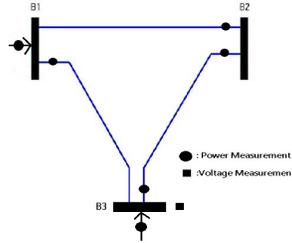

Figure 6: Configuration of a 3-bus system adopted from Example 5.5 in [1].

Fig. 6 shows the configuration of the 3-bus test system and its related measurements. Two bad data were implemented on the power flow measurement of the line $(3 - 2)$ and the power injection measurement on bus 3. Two scenarios have been considered: (i) multiple interacting, non-conforming bad data; (ii) multiple interacting, conforming bad data. In both scenarios random Gaussian noises have been added to the uncorrupted measurements. PSSE is carried out in each scenario by both the proposed joint optimization and traditional Gauss-Newton method with LNR test. The result is presented in Table **??**. As it can be seen, the proposed joint optimization algorithm is able to successfully identify the corrupted measurements even with the presence of interacting conforming bad data, while LNR method successfully identifies the bad data for nonconforming bad data but fails to detect any of them in the scenario when multiple interacting conforming bad data present.

Table 2: Effectiveness of proposed method in the presence of leverage points bad data

| Measu. | Bad data On inj1 | | | Bad data On Flow 1-2 | | |
|---|---|---|---|---|---|---|
| | $a$ | $o$ | $r^N$ | $a$ | $o$ | $r^N$ |
| Voltage 3 | $-0.0041$ | 0.0000 | **7.1958** | $-0.0023$ | 0.0000 | **7.5604** |
| Flow 1-2 | $-0.0001$ | 0.0000 | 1.6116 | **0.0598** | **0.0538** | 2.9419 |
| Flow 1-3 | 0.0003 | 0.0000 | 0.7135 | 0.0000 | 0.0000 | 0.9688 |
| Flow 3-1 | 0.0023 | 0.0000 | 0.5098 | 0.0029 | 0.0000 | 0.5772 |
| Flow 3-2 | 0.0042 | 0.0012 | 0.2900 | $-0.0021$ | 0.0000 | 1.261 |
| Flow 2-3 | $-0.0023$ | 0.0000 | 0.4884 | 0.0047 | 0.0003 | 1.3943 |
| Inj. 1 | **$-0.0351$** | **$-0.0306$** | 1.7143 | 0.0001 | 0.0000 | 3.1085 |
| Inj. 3 | 0.0001 | 0.0010 | 1.8129 | 0.0063 | 0.0002 | 1.3113 |

*4.4. Performance of the joint optimization algorithm against leverage point bad data*

This section investigates the effectiveness of the proposed joint optimization method in the presence of leverage point bad data. In order to do that, an 3-bus system example of leverage point bad data have been adapted from [36] and implemented on NESTA's 3-bus system to simulate leverage point bad data. In this example, the system has the same configuration as shown in Fig. 6 except that the power flow measurement on line 1-2 reads $P_{Flow_{12}}$ instead of $P_{Flow_{21}}$. In this example all the line assumed to have zero resistance and near 1 pu. reactance except for line 1-2 which is shortened to have a reactance of 0.1 pu. It has been shown in [36] that power injection measurement at bus 1 ($P_{inj_1}$) and the measurement flow on line 1-2 ($P_{Flow_{12}}$) become isolated leverage points. In order to investigate the effectiveness of joint optimization algorithm, two separate bad data have been added to the measurement flow 1-2 and power injection measurement 1. Table 2 shows the result of the test. The study clearly shows that joint optimization algorithm successfully identified leverage point bad data in both cases, while traditional LNR test fails to find either one.

*4.5. Sensitivity of the joint optimization algorithm to the Lagrangian multiplier*

Two case studies were carried out with $\lambda = 10^6$ and $\lambda = 10^{-6}$. In both case the objective function of the joint optimization compared to $\ell_1$-norm and $\ell_2$-norm of $z - h(x)$. The result has shown in Table **??**. It can be clearly seen that when $\lambda = 10^6$, joint optimization objective function reduces to WLSE while

Table 3: Effect of $\lambda$ on the bad data detection of the joint optimization algorithm

| Measu. | $a$ | $\lambda = 1.34$ | | $\lambda = 13.4$ | |
|---|---|---|---|---|---|
| | | $o$ | Detected | $o$ | Detected |
| $P_{\text{inj}_1}$ | **−0.0761** | **−0.0707** | **Yes** | 0.0000 | No |
| $P_{\text{inj}_2}$ | 0.0000 | 0.0000 | No | 0.0000 | No |
| $P_{\text{inj}_3}$ | 0.0000 | 0.0000 | No | 0.0000 | No |
| $Q_{\text{inj}_1}$ | **−0.0294** | **−0.0154** | **Yes** | 0.0000 | No |
| $Q_{\text{inj}_2}$ | 0.0000 | 0.0000 | No | 0.0000 | No |
| $Q_{\text{inj}_3}$ | 0.0000 | 0.0000 | No | 0.0000 | No |
| $P_{\text{flow}_{13}}$ | 0.0000 | 0.0000 | No | 0.0000 | No |
| $P_{\text{flow}_{32}}$ | 0.0000 | 0.0000 | No | 0.0000 | No |
| $P_{\text{flow}_{12}}$ | **0.0265** | **0.0180** | **Yes** | 0.0000 | No |
| $P_{\text{flow}_{31}}$ | 0.0000 | 0.0000 | No | 0.0000 | No |
| $P_{\text{flow}_{23}}$ | 0.0000 | 0.0000 | No | 0.0000 | No |
| $P_{\text{flow}_{21}}$ | **−0.0270** | **0.0181** | **Yes** | 0.0000 | No |
| $Q_{\text{flow}_{13}}$ | 0.0000 | 0.0000 | No | 0.0000 | No |
| $Q_{\text{flow}_{32}}$ | 0.0000 | 0.0000 | No | 0.0000 | No |
| $Q_{\text{flow}_{12}}$ | −0.0092 | 0.0000 | No | 0.0000 | No |
| $Q_{\text{flow}_{31}}$ | 0.0000 | 0.0000 | No | 0.0000 | No |
| $Q_{\text{flow}_{23}}$ | 0.0000 | 0.0000 | No | 0.0000 | No |
| $Q_{\text{flow}_{21}}$ | **0.0610** | **0.0570** | **Yes** | 0.0000 | No |

The bolded fonts indicate bad data.

when $\lambda = 10^{-6}$ it coincides to WLAV (at this case, $o = z - h(x)$.) Therefore, for a positive $\lambda$, the joint optimization will respond to Huber's M-estimator and optimal $\lambda$ can be chosen accordingly.

In the case study on NESTA 3-bus system, we have injected bad data on randomly selected real and reactive power injection and flow measurements. The error vector (including the noise) and the identified bad data for different $\lambda$ are shown in the Table 4. Note that there are total 21 measurements. The tables list only active and reactive power meters. It can be clearly seen in the table that $\lambda = 1.34$ is a better choice compared to $\lambda = 13.4$.

*4.6. More computing experiments*

In this section, performance of the proposed algorithm (joint state estimation and bad data detection) and the traditional WLSE solved by Gauss-Newton method with LNR detection test (Gauss-Newton + LNR). For the traditional method, an WLSE problem is first solved by Gauss-Newton. Then the LNR is identified. This measurement is excluded and the WLSE problem is solved again. The performance metrics here are the RMS of angle error (RMS-AE) and the RMS of voltage error (RMS-VE).

Multiple bad data are injected in test case systems. 10% of the randomly selected real and reactive power injection and flow measurements will be corrupted by bad data in addition to Gaussian noise. In all cases, the bad data are simulated by multiplying the true measurement by 1.2 and all the measurements are assumed to have noise.

Fig. 7 shows the voltage angle error of the two methods for NESTA's 14-bus system. It can be found that the joint optimization approach can give a much accurate estimation.

Table 4: Effect of $\lambda$ on the bad data detection of the joint optimization algorithm

| Measu. | $a$ | $\lambda = 1.34$ | | $\lambda = 13.4$ | |
|---|---|---|---|---|---|
| | | $o$ | Detected | $o$ | Detected |
| $P_{\text{inj}_1}$ | **−0.0761** | **−0.0707** | **Yes** | 0.0000 | No |
| $P_{\text{inj}_2}$ | 0.0000 | 0.0000 | No | 0.0000 | No |
| $P_{\text{inj}_3}$ | 0.0000 | 0.0000 | No | 0.0000 | No |
| $Q_{\text{inj}_1}$ | **−0.0294** | **−0.0154** | **Yes** | 0.0000 | No |
| $Q_{\text{inj}_2}$ | 0.0000 | 0.0000 | No | 0.0000 | No |
| $Q_{\text{inj}_3}$ | 0.0000 | 0.0000 | No | 0.0000 | No |
| $P_{\text{flow}_{13}}$ | 0.0000 | 0.0000 | No | 0.0000 | No |
| $P_{\text{flow}_{32}}$ | 0.0000 | 0.0000 | No | 0.0000 | No |
| $P_{\text{flow}_{12}}$ | **0.0265** | **0.0180** | **Yes** | 0.0000 | No |
| $P_{\text{flow}_{31}}$ | 0.0000 | 0.0000 | No | 0.0000 | No |
| $P_{\text{flow}_{23}}$ | 0.0000 | 0.0000 | No | 0.0000 | No |
| $P_{\text{flow}_{21}}$ | **−0.0270** | **0.0181** | **Yes** | 0.0000 | No |
| $Q_{\text{flow}_{13}}$ | 0.0000 | 0.0000 | No | 0.0000 | No |
| $Q_{\text{flow}_{32}}$ | 0.0000 | 0.0000 | No | 0.0000 | No |
| $Q_{\text{flow}_{12}}$ | −0.0092 | 0.0000 | No | 0.0000 | No |
| $Q_{\text{flow}_{31}}$ | 0.0000 | 0.0000 | No | 0.0000 | No |
| $Q_{\text{flow}_{23}}$ | 0.0000 | 0.0000 | No | 0.0000 | No |
| $Q_{\text{flow}_{21}}$ | **0.0610** | **0.0570** | **Yes** | 0.0000 | No |

The bolded fonts indicate bad data.

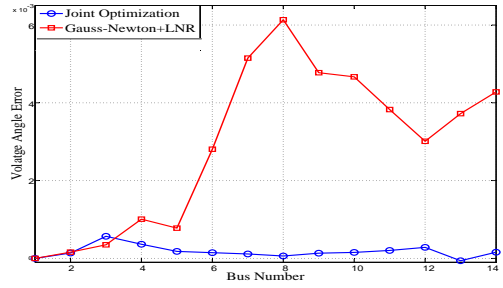

Figure 7: Voltage angle error on the IEEE-14 bus system in the presence of multiple bad data.

Table 5 represents the result of the comparison. The joint optimization adopts 3 iterations for SDP cuts generation. As a total, one SOCP relaxation problem and three more strengthened SOCP relaxation problems are computed. The estimation result shows better performance of the joint optimization method compared to the traditional LNR test. The computing time is also listed. The joint optimization method consumes more computing time for majority of the test cases. However, for the 1354-bus system, the Gauss-Newton based method consumed more than an hour while the joint optimization method consumed 28 minutes. For this case, the Gauss-Newton method does not converge and the maximum iteration number 100 is reached.

Table 5: Performance comparison between co-optimization algorithm and LNR test method in the presence of noise and multiple bad data

| Cases | Gauss-Newton+LNR | | | Joint Optimization | | |
|---|---|---|---|---|---|---|
| | RMS-AE | RMS-VE | Time (s) | RMS-AE | RMS-VE | Time (s) |
| nesta_case3_lmbd | 0.0060 | 0.0050 | 0.0382 | 0.0001 | 0.0009 | 0.3120 |
| nesta_case4_gs | 0.0005 | 0.0016 | 0.0211 | 0.0001 | 0.0014 | 0.3432 |
| nesta_case5_pjm | 0.0006 | 0.0002 | 0.0259 | 0.0022 | 0.0020 | 0.4836 |
| nesta_case6_ww | 0.0007 | 0.0046 | 0.2473 | 0.0004 | 0.0001 | 0.8424 |
| nesta_case14_ieee | 0.0048 | 0.0014 | 0.0665 | 0.0001 | 0.0001 | 1.5288 |
| nesta_case30_as | 0.0049 | 0.0002 | 0.1654 | 0.0002 | 0.0001 | 2.9640 |
| nesta_case30_fsr | 0.0085 | 0.0013 | 0.1861 | 0.0001 | 0.0004 | 2.9640 |
| nesta_case30_ieee | 0.0009 | 0.0029 | 0.1923 | 0.0001 | 0.0003 | 3.0888 |
| nesta_case57_ieee | 0.0130 | 0.0382 | 0.7765 | 0.0012 | 0.0004 | 6.0840 |
| nesta_case118_ieee | 0.1964 | 0.0301 | 3.4819 | 0.0055 | 0.0008 | 44.679 |
| nesta_case300_ieee | 0.2789 | 0.0893 | 22.963 | 0.0079 | 0.0025 | 370.84 |
| nesta_case1354_pegase | 0.2756 | 0.9103 | 3662.3 | 0.0119 | 0.0001 | 1689.1 |

## 5. Conclusion

A joint state estimation and bad data identification algorithm is introduced in this paper. The proposed algorithm uses the sparse signal characteristics to identify bad data. In order to be applied in AC network state estimation, our recently developed strengthened SOCP relaxation using LSE-based SDP cutting plane method is implemented to solve the joint state estimation and bad data identification problem. Numerical results from case studies demonstrate more accurate results in SOCP relaxation of state estimation, success of the algorithm for simultaneous state estimation and bad data identification and improved performance compared to largest normalized residual tests.